\documentclass{amsart}
\usepackage{amscd,amssymb,graphicx}
\author[Fialowski]{Alice Fialowski}
\address{
Alice Fialowski\\
E\"otv\"os Lor\'and University\\
Budapest, Hungary} \email{fialowsk@cs.elte.hu}
\author[Penkava]{Michael Penkava}
\address{
Michael Penkava\\
University of Wisconsin-Eau Claire\\
Eau Claire, WI 54702-4004} \email{penkavmr@uwec.edu}
\thanks{Research of the first author was partially supported by  OTKA grant K77757 and the second author by grants from the
University of Wisconsin-Eau Claire.}
\newtheorem{thm}{Theorem}[section]

\theoremstyle{definition}

\def \ph{\varphi}

\def\GL{\mbox{\bf GL}}

\def \diag{\operatorname {diag}}

\def \ra{\rightarrow}

\def \hom{\mbox{\rm Hom}}

\def \tns{\otimes}

\def \mcom{,\cdots,}
\def \k{\mbox{$\mathbb K$}}

\def \C{\mbox{$\mathbb C$}}
\def \Z{\mbox{$\mathbb Z$}}

\def\zt{\mbox{$\Z_2$}}

\def\ad{\operatorname{ad}}

\def\inv{^{-1}}

\def\im{\operatorname{Im}}

\def\p{\epsilon}

\def\ainf{\mbox{$A_\infty$}}
\def\linf{\mbox{$L_\infty$}}
\def\and{\mbox{ \rm and }}
\def\T{\mathcal T}
\def\TV{\T(V)}
\def\TW{\mbox{$\T(W)$}}

\def\pha#1#2{\ph^{#1}_{#2}}

\def\psa#1#2{\psi^{#1}_{#2}}

\def\inv{^{-1}}

\def\P{\mathbb P}

\def\xxx{d_{79}}
\def\yyy{d_{83}}
\input{xy}
\xyoption{all}
\begin{document}
\setlength{\multlinegap}{0pt}
\title[Nilpotent $4$-dimensional algebras]
{The Moduli space of $4$-dimensional nilpotent complex associative algebras}%

\address{}%
\email{}%

\thanks{}%

\date{\today}
%\dedicatory{}%
%\commby{}%
% ----------------------------------------------------------------
\begin{abstract}
In this paper, we study  $4$-dimensional nilpotent complex
associative algebras. This is a continuation of the study of the moduli space of
4-dimensional algebras. The non-nilpotent algebras were analyzed in an earlier paper.
Even though there are only 15 families of nilpotent 4-dimensional algebras, the complexity
of their behavior warranted a separate study, which we give here.
\end{abstract}
\maketitle
% ----------------------------------------------------------------

%\nocite{ps1,ps2,mp1,pv1}
%\input realnewintro.tex
\section{Construction of the algebras by extensions}
The authors and collaborators have been carrying out a construction of moduli spaces of low dimensional complex and real Lie and associative algebras in a series of papers.  Our method of constructing the moduli spaces of such algebras is based on the principal that algebras are either simple or can be constructed as extensions of lower dimensional algebras.  There is a classical theory of extensions, which was developed by many contributors going back as early as the 1930s.
In \cite{fp11}, we
gave a description of the theory of extensions of an algebra $W$ by an algebra $M$. Consider the diagram
$$
0\ra M\ra V\ra W\ra 0
$$
of associative \k-algebras, so that $V=M\oplus W$ as a \k-vector space, $M$ is an
ideal in the algebra $V$, and $W=V/M$ is the quotient algebra. Suppose that
$\delta\in C^2(W)=\hom(T^2(W),W)$ and $\mu\in C^2(M)$ represent the algebra structures on
$W$ and $M$ respectively. We can view $\mu$ and $\delta$ as elements of $C^2(V)$.
Let $T^{k,l}$ be the subspace of $T^{k+l}(V)$ given recursively
by
\begin{align*}
T^{0,0}&=\k\\
T^{k,l}&=M\tns T^{k-1,l}\oplus V\tns T^{k,l-1}.
\end{align*}
Let
$C^{k,l}=\hom(T^{k,l},M)\subseteq C^{k+l}(V)$.
If we denote the algebra structure on $V$ by $d$, we have
$$
d=\delta+\mu+\lambda+\psi,
$$
where $\lambda\in C^{1,1}$ and $\psi\in C^{0,2}$. Note that in this notation,
$\mu\in C^{2,0}$. Then the condition that $d$ is associative:  $[d,d]=0$ gives the
following relations:
\begin{align*}
[\delta,\lambda]+\tfrac 12[\lambda,\lambda]+[\mu,\psi]&=0,
\quad\text{The Maurer-Cartan equation}\\
[\mu,\lambda]&=0,\quad\text{The compatibility condition}\\
[\delta+\lambda,\psi]&=0,\quad\text{The cocycle condition}
\end{align*}
Since $\mu$ is an algebra structure, $[\mu,\mu]=0$, so if we define
$D_\mu$ by $D_\mu(\ph)=[\mu,\ph]$, then $D^2_\mu=0$.
Thus $D_\mu$ is a differential on $C(V)$.
Moreover $D_\mu:C^{k,l}\ra C^{k+1,l}$. Let
\begin{align*}
Z_\mu^{k,l}&=\ker(D_\mu:C^{k,l}\ra C^{k+1,l}),\quad\text{the $(k,l)$-cocycles}\\
B_\mu^{k,l}&=\im(D_\mu:C^{k-1,l}\ra C^{k,l}),\quad\text{the $(k,l)$-coboundaries}\\
H_\mu^{k,l}&=Z_\mu^{k,l}/B_\mu^{k,l},\quad\text{the $D_\mu$-cohomology}
\end{align*}

Then the compatibility condition means that $\lambda\in Z^{1,1}$.
If we define $D_{\delta+\lambda}(\ph)=[\delta+\lambda,\ph]$, then it is not
true that $D^2_{\delta+\lambda}=0$, but
$D_{\delta+\lambda}D_\mu=-D_{\mu}D_{\delta+\lambda}$, so that $D_{\delta+\lambda}$ descends
to a map $D_{\delta+\lambda}:H^{k,l}_\mu\ra H^{k,l+1}_\mu$, whose square is zero, giving
rise to the $D_{\delta+\lambda}$-cohomology $H^{k,l}_{\mu,\delta+\lambda}$.
If the pair $(\lambda,\psi)$ gives rise to an algebra $d$, and $(\lambda,\psi')$
give rise to another algebra $d'$, then if we express $\psi'=\psi+\tau$, it is
easy to see that $[\mu,\tau]=0$, and $[\delta+\lambda,\tau]=0$, so that the image $\bar\tau$
of $\tau$ in $H^{0,2}_\mu$ is a $D_{\delta+\lambda}$-cocycle, and thus $\tau$ determines
an element $\{\bar\tau\}\in H^{0,2}_{\mu,\delta+\lambda}$.

If $\beta\in C^{0,1}$, then $g=\exp(\beta):\TV\ra\TV$ is given by
$g(m,w)=(m+\beta(w),w)$. Furthermore, $g^*=\exp(-\ad_{\beta}):C(V)\ra C(V)$ satisfies
$g^*(d)=d'$, where $d'=\delta+\mu+\lambda'+\psi'$ with
$\lambda'=\lambda+[\mu,\beta]$ and $\psi'=\psi+[\delta+\lambda+\tfrac12[\mu,\beta],\beta]$.
In this case, we say that $d$ and $d'$ are equivalent extensions in the restricted sense.
Such equivalent extensions are also equivalent as algebras on $\TV$.
Note that $\lambda$
and $\lambda'$ differ by a $D_\mu$-coboundary, so $\bar\lambda=\bar\lambda'$ in
$H^{1,1}_\mu$. If $\lambda$ satisfies the MC-equation for some $\psi$, then
any element $\lambda'$ in $\bar\lambda$ also gives a solution of the MC equation,
for the $\psi'$ given above. The cohomology classes of those $\lambda$ for which
a solution of the MC equation exists, determine distinct restricted equivalence classes
of extensions.

Let $G_{M,W}=\GL(M)\times\GL(W)\subseteq\GL(V)$. If $g\in G_{M,W}$ then $g^*:C^{k,l}\ra
C^{k,l}$, and $g^*:C^k(W)\ra C^k(W)$, so $\delta'=g^*(\delta)$ and $\mu'=g^*(\mu)$
are algebra structures on $\T(M)$ and $\TW$ respectively.
The group $G_{\delta,\mu}$ is the
subgroup of $G_{M,W}$ consisting of those elements $g$ such that $g^*(\delta)=\delta$
and $g^*(\mu)=\mu$. Then $G_{\delta,\mu}$ acts on
the restricted equivalence classes of extensions, giving the equivalence classes
of general extensions. Also $G_{\delta,\mu}$
acts on $H^{k,l}_\mu$, and induces an action on the classes $\bar\lambda$ of $\lambda$
giving a solution to the MC equation.

Next, consider the group $G_{\delta,\mu,\lambda}$ consisting
of the automorphisms $h$ of $V$ of the form $h=g\exp(\beta)$, where
$g\in G_{\delta,\mu}$, $\beta\in C^{0,1}$ and $\lambda=g^*(\lambda)+[\mu,\beta]$.
If $d=\delta+\mu+\lambda+\psi+\tau$, then $h^*(d)=\delta+\mu+\lambda+\psi+\tau'$ where
\begin{equation*}
\tau'=g^*(\psi)-\psi+[\delta+\lambda-\tfrac12[\mu,\beta],\beta]+g^*(\tau).
\end{equation*}
Thus the group $G_{\delta,\mu,\lambda}$ induces an action on $H^{0,2}_{\mu,\delta+\lambda}$
given by $\{\bar\tau\}\ra\{\overline{\tau'}\}$.

The general group of equivalences of extensions of the algebra structure $\delta$ on $W$
by the algebra structure $\mu$ on $M$ is given by the group of automorphisms of $V$ of
the form $h=\exp(\beta)g$, where $\beta\in C^{0,1}$ and $g\in G_{\delta,\mu}$. Then we
have the following classification of such extensions up to equivalence.
\begin{thm}[\cite{fp11}]
The equivalence classes of extensions of $\delta$ on $W$ by $\mu$ on $M$ is classified
by the following:
\begin{enumerate}
\item Equivalence classes of $\bar\lambda\in H^{1,1}_\mu$ which satisfy the MC equation
\begin{equation*}
[\delta,\lambda]+\tfrac12[\lambda,\lambda]+[\mu,\psi]=0
\end{equation*}
for some $\psi\in C^{0,2}$, under the action of the group $G_{\delta,\mu}$.
\item Equivalence classes of $\{\bar\tau\}\in H^{0,2}_{\mu,\delta+\lambda}$ under the
action of the group $G_{\delta,\mu,\lambda}$.
\end{enumerate}
\end{thm}
Equivalent extensions will give rise to equivalent algebras on $V$, but it may
happen that two algebraic structures arising from nonequivalent extensions are equivalent.
This is because the group of equivalences of extensions is the group of invertible
block upper
triangular matrices on the space $V=M\oplus W$, whereas the equivalence
classes of algebraic structures on $V$ are given by the group of all invertible
matrices, which is larger.

The fundamental theorem of finite dimensional algebras allows us to restrict our
consideration of extensions to two cases. First, we can consider those extensions
where $\delta$ is a semisimple algebra structure on $W$, and $\mu$ is a nilpotent
algebra structure on $M$. In this paper let $\k=\C$. Then we can
also assume that $\psi=\tau=0$. Thus the classification of the extension reduces
to considering equivalence classes of $\lambda$.

Secondly, we can consider extensions
of the trivial algebra structure $\delta=0$ on a 1-dimensional space $W$ by
a nilpotent algebra $\mu$. This
is because a nilpotent algebra has a codimension 1 ideal $M$, and the restriction
of the algebra structure to $M$ is nilpotent. However, in this case, we cannot assume
that $\psi$ or $\tau$ vanish,
so we need to use the classification theorem above to determine the
equivalence classes of extensions. In many cases, in solving the MC equation for
a particular $\lambda$, if there is any $\psi$ yielding a solution, then $\psi=0$
also gives a solution. So the action of $G_{\delta,\mu,\lambda}$ on $H^{0,2}_\mu$
takes on a simpler form than the general action we described above.

In addition to the complexity which arises because we cannot take the cocycle term
$\psi$ in the extension to be zero, there is another issue that complicates the construction
of the extensions. If an algebra is not nilpotent, then it has a maximal nilpotent ideal
which is unique, and it can be constructed as an extension of a semisimple algebra by this
unique ideal. Both the semisimple and nilpotent parts in this construction are completely
determined by the algebra. Therefore, a classification of extensions up to equivalence of
extensions will be sufficient to classify the algebras. This means that the equivalence classes
of the module structure $\lambda$ determine the algebras up to isomorphism.

For nilpotent algebras, we don't have this assurance. The same algebra structure may arise by
extensions of the trivial algebra structure on a 1-dimensional space by two different nilpotent
algebra structures on the same $n-1$-dimensional space.

In addition, the deformation theory of nilpotent algebras is far more involved than the deformation
theory of the non nilpotent algebras. 

In this paper, we study the complex 4-dimensional nilpotent algebras.
In \cite{degraaf}, the following idea for construction of nilpotent algebras was discussed, and since it is simpler
than  the general construction, we outline the idea here.

First, if a nonzero algebra is nilpotent, then it has a nontrivial ideal $M$ which has the property that the product of any
element in $M$ with an arbitrary element is zero. We call such an ideal completely trivial.  There is a unique ideal which
is maximal in the set of completely trivial ideals, which we call the
\emph{kernel} of the algebra.  If we call the quotient of the algebra
by its kernel the \emph{core} of the algebra, then we see that any nilpotent algebra determines a unique kernel and core, and the algebra
is given by an extension of its core by the kernel in a particularly simple manner.  In the language we introduced above,
we have $\lambda=0$ and $\mu=0$, so that the compatibility relation and
the Maurer Cartan equation are satisfied trivially.
In terms of the coboundary operator $D_\delta$ given by
$D_\delta(\ph)=0$, we get that the cocycle $\psi$ is actually a $D_\delta$-cocycle,
and equivalent cocycles determine equivalent algebras, so that the algebras are classified by the action of the group $G_{\delta,\mu,\lambda}$ on
the $D_\delta$-cohomology classes.

All of this holds for any completely trivial ideal, and it is more convenient to study the case when we don't assume $M$ is the kernel, even though
this means that some of the algebras will be constructed in a non unique manner.  In particular, let us consider the case where $\delta$ vanishes
as well as $\mu$ and $\lambda$.  Then the group $G_{\mu,\delta,\lambda}$ coincides with the group $G_{M,W}$.  If we express $M=\langle e_1\rangle$,
$W=\langle e_2,\dots,e_n\rangle$, then we can express
$\psi=\psa{i+1,j+1}1c_{i,j}$, and an element $g$ in $G_{M,W}$ is given in block form by a matrix
$G=\left[\begin{smallmatrix}g_{1,1}&0\\0&G_2\end{smallmatrix}\right]$. If $C=\left[c_{i,j}\right]$ is the $n\times n$ matrix determined by the coefficients of $\psi$,
and $\psi'=G^*(\psi)$ is the cocycle determined by the action of $G$ on $\psi$, and $C'=\left[c'_{i,j}\right]$ is the corresponding matrix of coefficients
of $\psi'$, then $C'=g\inv_{1,1}G^TCG$, so that $C'$ and $C$ are \emph{cogredient} matrices.  Therefore, the classification of cogredient matrices,
or equivalently, complex bilinear forms, is a useful component of the construction.

In the next section, we classify complex bilinear forms on a 2-dimensional complex vector space and then in the following section, we use this
classification to give a construction of the nilpotent 3-dimensional complex algebras, which agrees with the classification in \cite{fpp1}.
Then in the succeeding sections we use the same idea to construct the nilpotent 4-dimensional complex algebras.

\section{Classification of bilinear forms on a 2-dimensional complex vector space.}
A complete classification of complex bilinear forms on a finite
dimensional space was given e.g.
in \cite{ho-se1}. However, this classification does not
involve a stratification of the moduli space by complex projective orbifolds, which we believe is the right
way to understand this classification. Certainly, it is the correct point of view for the purpose for which we will apply it.
In this paper, we will give a stratification of the moduli spaces of 2 and 3 dimensional complex bilinear forms by projective
orbifolds.

Let $\beta$ be a bilinear form on $\C^2$ which is given by the matrix $B=\left(b_{ij}\right)$ where
$b_{ji}=\beta(e_i,e_j)$ in terms of some basis $\langle e_1,e_2\rangle$. We say that two matrices
$B$ and $C$ are \emph{cogredient} if there is a nonsingular matrix $P$ such that $C=P^TBP$,
This is precisely the condition that $B$ and $C$ represent $\beta$ with respect to different bases.

Then every nontrivial bilinear form is given up to equivalence by either a matrix of the form
$B(p: q)=\left[\begin{matrix}1&p\\q&0\end{matrix}\right]$ or $C=\left[\begin{matrix}0&1\\-1&0\end{matrix}\right]$.
Moreover, the matrices of the form $B(p:q)$ form a projective family parameterized by the orbifold $\P^1/\Sigma_2$,
where the action of the symmetric group $\Sigma_2$ on $\P^1$ is given by interchanging the projective coordinates $(p: q)$,
so that $B(p: q)\sim B(q: p)$ where $\sim$ stands for the equivalence of cogredient matrices.

To see why this is true, we note first that if $\beta$ is nontrivial and $\beta(u,u)=0$ for all $u\in\C^2$, then
$\beta(v,u)=-\beta(u,v)$ for all $u,v\in\C^2$, since $0=\beta(u+v,u+v)=\beta(u,v)+\beta(v,u)$. Since $\beta$ is nontrivial,
there must be some $u$ and $v$ such that $\beta(v,u)=1$, and the matrix of $\beta$ with respect to this basis is $C$.

Thus we may assume there is some $u$ in $\C^2$ such that $\beta(u,u)=1$. We claim that there is a nonzero
vector $v$ such that $\beta(v,v)=0$. For suppose that we choose any vector $v$ which is linearly independent of $u$.
If $\beta(v,v)\ne 0$, let $w=u+xv$, and then we compute
$$
\beta(u,w)=1+(\beta(u,v)+\beta(v,u))x+\beta(v,v)x^2.
$$
This is a quadratic equation in $x$ which has a nontrivial solution $x\ne 0$.  In terms of the basis $u,w$, the matrix
of $\beta$ has the form $B(p: q)$ for some $(p: q)$. Moreover, it is easy to check that $B(p: q)\sim B(x: y)$ precisely when
$(p: q)=(x: y)$ or $(p: q)=(y: x)$.

If one studies the miniversal deformation of the moduli space of bilinear forms, one discovers that elements of the form
$B(p: q)$ have smooth deformations in a neighborhood of the element $B(p: q)$ and do not have jump deformations, while the element
$C$ has a jump deformation to $B(1: -1)$ and smooth deformations in a
neighborhood of $B(1: -1)$. We shall see that the deformations of the moduli
space of bilinear forms are reproduced in the deformations of the moduli space of three dimensional complex algebras which
are determined by these bilinear forms. For basic notions of
deformations see \cite{fp10}.

\section{Nilpotent 3-dimensional complex associative algebras}
The moduli space of complex 3-dimensional associative algebras was constructed in \cite{fpp1}. Here we wish to construct the nilpotent algebras
using the following observations.  In \cite{degraaf}, the author
discusses what he calls central extensions of an associative algebra.
These are
extensions such that the ideal in the extended algebra has trivial multiplication with the entire algebra. This language derives from the similar
construction of central extensions of Lie algebras. In this section, we
show how to use the classification of bilinear forms as a tool in constructing the 3-dimensional nilpotent algebras.

There is only one nontrivial nilpotent 2-dimensional complex algebra,
represented by $\delta=\psa{33}2$, on the space $W=\langle e_2,e_3\rangle$.
If we extend $\delta$ by a completely trivial ideal $M=\langle e_1\rangle$, then this extension is given by a cocycle
$\psi=\psa{22}1c_1+\psa{23}1c_2+\psa{32}1c_3+\psa{33}1c_4$. The cocycle condition $[\delta,\psi]=0$ gives $c_1=0$ and $c_2=c_3$. Moreover,
if $\beta=\pha21b_1+\pha31b_2\in C^{1,0}$, then $[\delta,\beta]=-\psa{33}1b_1$, which means that up to a coboundary term, we can assume
$\psi=(\psa{2,3}1+\psa{3,2}1)c_3$.  Applying an element of the group $G_{\delta,\mu,\lambda}$ to $\psi$ replaces $\psi$ by an arbitrary
nonzero multiple, so this means we only have two cases to study, when
$c_3=1$ or $c_3=0$.  The first case gives the algebra
$d=\psa{33}2+\psa{23}1+\psa{32}1$, which is equivalent to the algebra $d_{19}$ on our list of 3-dimensional algebras.  The second case
gives $d=\psa{33}2$, which is equivalent to $d_{20}(0: 0)$ on our our list, but note that this algebra has a kernel which is 2-dimensional,
so this algebra would arise in a different fashion as well.

Next, consider the trivial nilpotent 2-dimensional complex algebra given by $\delta=0$. In this case the cocycle condition on
$\psi=\psa{22}1c_1+\psa{23}1c_2+\psa{32}1c_3+\psa{33}1c_4$ is trivial.  Let $C=\left[\begin{smallmatrix}c_1&c_2\\c_3&c_4\end{smallmatrix}\right]$.
An element $g$ of the group $G_{\delta,\mu,\lambda}$ is given by an arbitrary invertible matrix $G$ of the block diagonal form
$G=\diag(g_1,G_2)$ where $g_1\ne 0$ and $G_2$ is an invertible $2\times2$ matrix.  The action of $g$ on $\psi$ transforms $\psi$ into
the element $\psi'$ whose matrix is $g\inv_1G_2^{T}CG_2$, which means that if the matrices representing $\psi$ and $\psi'$ are cogredient
then they determine equivalent algebras. Using our classification of cogredient matrices, we find that we obtain a family
$d_{20}(p: q)$ given by the matrix $B(p: q)$ in our classification of bilinear forms on a 2-dimensional
complex vector space, and the algebra $d_{21}$, given by the matrix $C$ in this classification.

This gives all of the nontrivial 3-dimensional  nilpotent complex associative algebras in a very simple fashion.

\section{Classification of bilinear forms on a 3-dimensional complex vector space}
If $\beta$ is a bilinear form on an $n$-dimensional space $U$, then we say that $\beta$ is decomposable if $U=V\oplus W$, where $V$ and $W$ are nontrivial
subspaces of $U$ satisfying $\beta(V,W)=\beta(W,V)=0$.  Using the classification of 2-dimensional complex bilinear forms, it is easy to see that
every matrix $B$ representing a nontrivial  bilinear form on a complex 3-dimensional space is cogredient to a matrix of one of the six types given below.
\begin{align*}
B_1(p:q)=\left[\begin{matrix}1&0&0\\0&1&p\\0&q&0\end{matrix}\right],\qquad
B_2(p:q)=\left[\begin{matrix}0&0&0\\0&1&p\\0&q&0\end{matrix}\right]\\
B_3=\left[\begin{matrix}1&0&0\\0&0&1\\0&-1&0\end{matrix}\right],\qquad
B_4=\left[\begin{matrix}0&0&0\\0&0&1\\0&-1&0\end{matrix}\right]\\
B_5=
\left[\begin{matrix}0&1&0\\0&0&1\\0&-1&0\end{matrix}\right],\qquad
B_6=
\left[\begin{matrix}1&1&0\\0&1&1\\0&1&0\end{matrix}\right].
\end{align*}
The first four matrices correspond to the decomposable bilinear forms, and it follows from our classification of bilinear forms on $\C^2$ that the matrix of every decomposable bilinear form is cogredient to one of these six
matrix types. It is not as easy to see that the matrix of any indecomposable bilinear form is cogredient to one of the last two matrices. In fact, from the 
results of \cite{ho-se1}, it follows that when $n$ is odd, there are
exactly two nonequivalent $n\times n$ indecomposable matrices.

The matrix $B_1(0: 0)$ is cogredient to the matrix $B_3(1: 1)$. Moreover, the families $B_1(p: q)$ and $B_3(p:q)$ are parametrized by
$\P^1/\Sigma_2$, so they determine projective orbifolds.  Other than these identifications, all of the matrices represent distinct equivalence classes.
As a consequence, we see that the moduli space of complex bilinear forms on a 3-dimensional vector space is stratified by projective orbifolds. It is also
true that the deformations of the elements in the moduli space are either given by jumps between the strata, by smooth deformations which factor through
a jump deformation, or by smooth deformations along a stratum.  This pattern is consistent with the patterns that we have observed in moduli spaces of
algebras.  

The decomposition of the moduli space of $3\times 3$ matrices we have given above has the advantages that it gives a stratification by projective orbifolds and
moreover, it is easy to see which matrices are indecomposable.
Nevertheless, even this decomposition is not the correct one from the
point of view of deformation theory. Consider the following 6 matrices:
\begin{align*}
C_1(p:q)&=\left[\begin{matrix}0&0&q\\0&1&1\\p&0&1\end{matrix}\right],\qquad
C_2(p:q)=\left[\begin{matrix}0&0&q\\0&0&0\\p&0&1\end{matrix}\right]\\
C_3&=\left[\begin{matrix}0&1&0\\-1&0&0\\0&0&1\end{matrix}\right],\qquad
C_4=\left[\begin{matrix}0&1&1\\-1&0&0\\0&0&0\end{matrix}\right]\\
C_5&=\left[\begin{matrix}0&1&0\\1&1&0\\0&0&1\end{matrix}\right],\qquad
C_6=\left[\begin{matrix}0&1&0\\-1&0&0\\0&0&0\end{matrix}\right].
\end{align*}
The matrices above correspond to the algebras $d_{78}(p:q)$, $d_{86}(p:q)$, $d_{81}$, $d_{84}$, $d_{85}$ and $d_{87}$, which we will describe below.
In most cases, there is a one to one correspondence between the $B$ matrices above and the $C$ matrices. However, there are a couple of very important 
exceptions.  The matrix $B_1(p:q)$ is cogredient to $C_1(p:q)$ except for $(1:1)$ and $(0:0)$. The matrix $B_1(1:1)$ is cogredient to $C_5$, while the
matrix $B_1(0:0)$ is cogredient to $C_2(1:1)$. The matrix $B_2(p:q)$ is always cogredient to $C_2(p:q)$, $B_3$ is cogredient to $C_3$, $B_4$ is cogredient
to $C_6$, $B_5$ is cogredient to $C_5$, and $B_6$ is cogredient to $C_1(1:1)$. 

The interchange of $B_1(1:1)$ with $B_6$ seems a bit strange because $B_6$ is indecomposable, so it would seem that the family of decomposables $B_1(p:q)$ is more natural than the latter stratification.  However, it turns out that decomposability/indecomposability is not preserved under deformations. The stratification given by the $C$ matrices is consistent with deformation theory in that deformations occur along a stratum, jump deformations to elements of different strata, and smooth deformations along another stratum which ``factor through a jump deformation''.  

\section{Stratification of the nilpotent algebras}

The nilpotent algebras are divided into 15 different strata. There are 4 projective 1-parameter  families of algebras: $d_{75}(p:q)$, $d_{78}(p:q)$, $\yyy(p:q)$ and $d_{86}(p:q)$, where $(p:q)$ is a projective coordinate. For some of these families, there is also an action of the symmetric
group $\Sigma_2$, given by permutation of the coordinates, so that the algebra associated to the parameter $(p:q)$ is isomorphic to the algebra
given by the parameter $(q:p)$. In this case, the family is parameterized by $\P^1/\Sigma_2$. Otherwise the family is parameterized by $\P^1$.
The 11 algebras $d_{73}$, $d_{74}$, $d_{76}$, $d_{77}$, $d_{80}$, $d_{81}$, $d_{82}$, $\xxx$, $d_{84}$, $d_{85}$ and $d_{87}$ each determine
a stratum consisting of a single algebra. We will give a description of each of the strata below.

\begin{table}[h,t]
\begin{center}
\begin{tabular}{lcccc}
Algebra &$H^0$&$H^1$&$H^2$&$H^3$\\ \hline \\
$d_{73}=\psa{31}2+\psa{13}2+\psa{33}1+\psa{34}2+\psa{44}2$&$0|2$&$2|0$&$0|3$&$4|0$\\
$d_{74}=\psa{33}2+\psa{41}2+\psa{43}1+\psa{14}2+\psa{34}1+\psa{44}3$&$0|4$&$4|0$&$0|4$&$4|0$\\
$d_{75}(p:q)=\psa{33}2+p\psa{43}1+q\psa{34}1+\psa{44}2$&$0|2$&$3|0$&$0|5$&$8|0$\\
$d_{75}(1:1)=\psa{33}2+\psa{34}2+\psa{44}1$&$0|2$&$3|0$&$0|6$&$10|0$\\
$d_{75}(1:-1)=\psa{42}1+\psa{24}3+\psa{44}1$&$0|2$&$3|0$&$0|6$&$11|0$\\
$d_{75}(1:0)=\psa{33}2+\psa{43}1+\psa{44}2$&$0|2$&$3|0$&$0|5$&$8|0$\\
$d_{75}(0:0)=\psa{33}2+\psa{44}2$&$0|4$&$8|0$&$0|17$&$41|0$\\
$d_{76}=\psa{31}2+\psa{13}2+\psa{33}1+\psa{44}2$&$0|4$&$5|0$&$0|6$&$8|0$\\
$d_{77}=\psa{31}2+\psa{13}2+\psa{33}1+\psa{34}2$&$0|2$&$3|0$&$0|6$&$10|0$\\
$d_{78}(p:q)=\psa{33}2+p\psa{41}2+q\psa{14}2+\psa{34}2+\psa{44}2$&$0|2$&$3|0$&$0|6$&$12|0$\\
$d_{78}(1:1)=\psa{33}2+\psa{41}2+\psa{14}2+\psa{34}2+\psa{44}2$&$0|2$&$3|0$&$0|6$&$12|0$\\
$d_{78}(1:-1)=\psa{33}2+\psa{41}2-\psa{14}2+\psa{34}2+\psa{44}2$&$0|2$&$3|0$&$0|6$&$12|0$\\
$d_{78}(1:0)=\psa{33}2+\psa{41}2+\psa{34}2+\psa{44}2$&$0|2$&$3|0$&$0|8$&$23|0$\\
$d_{78}(0:0)=\psa{33}2+\psa{34}2+\psa{44}2$&$0|2$&$6|0$&$0|15$&$37|0$\\
$\xxx=\psa{31}2+\psa{13}2+\psa{33}1$&$0|4$&$6|0$&$0|10$&$18|0$\\
$d_{80}=\psa{33}2+\psa{43}1-\psa{34}1+\psa{44}2$&$0|2$&$4|0$&$0|9$&$17|0$\\
$d_{81}=-\psa{31}2+\psa{13}2+\psa{44}2$&$0|2$&$5|0$&$0|10$&$18|0$\\
$d_{82}=\psa{33}2+\psa{43}1+\psa{34}1+\psa{44}2$&$0|4$&$6|0$&$0|10$&$18|0$\\
$\yyy(p:q)=p\psa{42}1+q\psa{24}1+\psa{44}1+\psa{44}3$&$0|2$&$5|0$&$0|9$&$16|0$\\
$\yyy(1:1)=\psa{42}1+\psa{24}1+\psa{44}1+\psa{44}3$&$0|4$&$7|0$&$0|12$&$21|0$\\
$\yyy(1:-1)=\psa{42}1-\psa{24}1+\psa{44}1+\psa{44}3$&$0|2$&$5|0$&$0|10$&$18|0$\\
$\yyy(1:0)=\psa{42}1+\psa{44}1+\psa{44}3$&$0|2$&$5|0$&$0|11$&$21|0$\\
$\yyy(0:0)=\psa{44}1+\psa{44}3$&$0|4$&$10|0$&$0|28$&$82|0$\\
$d_{84}=-\psa{31}2+\psa{13}2+\psa{14}2$&$0|2$&$4|0$&$0|10$&$26|0$\\
$d_{85}=\psa{31}2+\psa{13}2+\psa{33}2+\psa{44}2$&$0|4$&$7|0$&$0|14$&$28|0$\\
$d_{86}(p:q)=p\psa{42}1+q\psa{24}1+\psa{44}1$&$0|2$&$6|0$&$0|15$&$37|0$\\
$d_{86}(1:1)=\psa{42}1+\psa{24}1+\psa{44}1$&$0|4$&$8|0$&$0|17$&$41|0$\\
$d_{86}(1:-1)=\psa{42}1-\psa{24}1+\psa{44}1$&$0|2$&$6|0$&$0|15$&$38|0$\\
$d_{86}(1:0)=\psa{42}1+\psa{44}1$&$0|2$&$6|0$&$0|19$&$52|0$\\
$d_{86}(0:0)=\psa{44}1$&$0|4$&$10|0$&$0|28$&$82|0$\\
$d_{87}=-\psa{31}2+\psa{13}2$&$0|2$&$8|0$&$0|17$&$42|0$\\
\end{tabular}
\end{center}
\label{d73-d87}
\caption{The cohomology of the algebras $d_{73}\dots d_{87}$}
\end{table}

\section{The fifteen families of algebras}
\subsection{Type 73}
The algebra $d_{73}=\psa{31}2+\psa{13}2+\psa{33}1+\psa{34}2+\psa{44}2$ has a 1-dimensional kernel $M=\langle e_2\rangle$.
Its core is the 3-dimensional nilpotent algebra $d_{20}(0:0)=\psa{33}2$. The algebra $d_{73}$ is isomorphic to its opposite algebra.
We compute
\begin{equation*}
H^2(d_{73})=\langle \delta^1,\delta^2,\delta^3\rangle,
\end{equation*}
where
\begin{align*}
\delta^1&=-\psa{1,1}2-\psa{2,3}2-\psa{3,2}2+\psa{3,4}1+\psa{4,4}1\\
\delta^2&=-\psa{1,1}2+\psa{2,3}2+\psa{3,2}2+2\psa{3,3}3+\psa{3,4}4+\psa{4,3}4\\
\delta^3&=-2\psa{1,3}4-2\psa{3,1}4+2\psa{3,2}2+\psa{1,4}1+3\psa{2,4}2+\psa{3,4}3-2\psa{3,4}4+\psa{4,1}1+3\psa{4,2}2+\psa{4,3}3.
\end{align*}
More precisely, the cohomology classes of these algebras give a basis of $H^2$, but it is convenient to identify the cohomology
classes with \emph{representative cocycles}, which give a pre-basis of the cohomology.

The third order deformation is versal, with four relations
\begin{equation*}
t_3^2=0,\quad t_2t_3=0,\quad t_2t_3=0,\quad t_2^2(t_1+t_3)=0.
\end{equation*}
Note that the number of relations is precisely the dimension of $H^3$, and that the third relation is redundant.  From the first relation, we see that in any deformation, $t_3=0$. Since $t_3$ vanishes by the first relation, the second and third relations give no additional constraints. Finally, the fourth relation gives that either $t_1=0$ or $t_2=0$. When $t_1=0$, we obtain a 1-parameter deformation
\begin{align*}
d_{t_2}=&\psa{31}2+\psa{13}2+\psa{33}1+\psa{34}2+\psa{44}2-t_2\psa{11}2+t_2\psa{23}2
+t_2\psa{32}2+2t_2\psa{33}3+t_2\psa{34}4\\&+t_2\psa{43}4-t_2^2\psa{11}1-t_2^2\psa{12}2
-t_2^2\psa{21}2-t_2^2\psa{13}3-t_2^2\psa{31}3-t_2^2\psa{14}4-t_2^2\psa{41}4
\end{align*}
which is a deformation to a neighborhood of the algebra $d_{37}(1:1)$. In fact, one can show that $d_{t_2}\sim d_{37}(x:y)$ where $t_2=-\frac{(x-y)^2}{xy}$, whenever $t_2\ne 0$. From this expression, we see that as $t_2\ra 0$, the algebra $d_{t_2}(x:y)\ra d_{37}(1: 1)$. However, when $t_2=0$, we obtain the original algebra $d_{73}$. What happens is that the matrix which expresses the isomorphism between the algebras $d_{t_2}$ and $d_{37}(x:y)$ becomes singular when
$x=y$.

The second solution of the relations, given by setting $t_2=0$, gives a
1-parameter deformation in a neighborhood of $d_{38}(1:1)$. As in the case of the first
solution, we do not obtain a deformation to the algebra $d_{38}(1:1)$.

There are no
jump deformations to either $d_{37}(1:1)$ or $d_{38}(1:1)$, which we found surprising, since in our previous constructions, we have always
observed that when a member of one stratum deforms smoothly to a
neighboneighborhood of a point in a different stratum, then the smooth deformation factored
through a jump deformation to that point, which just means that there is a jump deformation to the point.

\subsection{Type 74}
The algebra $d_{74}=\psa{33}2+\psa{41}2+\psa{43}1+\psa{34}1+\psa{14}2+\psa{44}3$ has a 1-dimensional kernel $M=\langle e_2\rangle$.
Its core is the 3-dimensional nilpotent algebra $d_{19}=\psa{43}1+\psa{34}1+\psa{44}3$. This algebra is commutative.

When we compute the miniversal deformation, we obtain that the fourth order deformation is versal,
and all the relations vanish. The miniversal deformation gives jump deformations to
$d_2$, $d_6$, $d_7$, $d_8$, $d_9$, $d_{31}$, $d_{32}$, $d_{35}$, $d_{36}$, $d_{49}$, and $d_{50}$.

It is interesting to note that
all the algebras to which $d_{74}$ deforms are also commutative. It is well-known that in order for an algebra to deform to a commutative
algebra, it must be commutative. However, it is possible, and even common, for commutative algebras to deform to noncommutative algebras.

\subsection{Type $\mathbf{75(p:q)}$}
The family of algebras $d_{75}(p:q)$ are generically given by the
algebras
$\psa{33}2+p\psa{43}1+q\psa{34}1+\psa{44}2$.  However, for certain special values of the parameters,
the algebra is given by a different algebraic structure.  We have $d_{75}(1: -1)=\psa{43}1+\psa{34}2+\psa{44}1$ and
$d_{75}(1: 1)=\psa{33}2+\psa{34}2+\psa{44}1$.  The family is parameterized by $\P^1/\Sigma_2$, so that $d_{75}(p:q)\sim d_{75}(q:p)$.
It might be possible to find a family  $d(p:q)$ of algebras which give a parametrization of the entire family, but we did not
discover such a representation.  In all cases, the algebra $d_{75}(p:q)$ has a 2-dimensional kernel $M=\langle e_1,e_2\rangle$ and
its core is the trivial nilpotent 2-dimensional algebra.

For projective families of algebras, generically, the dimension of $H^n$ does not vary, and there is a generic pattern for the deformations.  There are a few
values of $(p:q)$ for which the pattern varies from the generic case, in that the dimension of $H^n$ may be larger, and there may be extra deformations.
Usually, these special values are $(1:1)$, $(1:-1)$, $(1:0)$ and $(0:0)$.  In fact, $(0:0)$ is always special, and there are always jump deformations from
the algebra correcponding to $(0:0)$, to every member in the family.  When there is no
action of $\Sigma_2$ on the family, so the aalgebra corresponding
to $(p:q)$ is not generically isomorphic to that for $(q:p)$, the value $(0:1)$ is also special.  Occasionally, there are additional special values.
\subsubsection{Generic Case}
In the case of $d_{75}(p:q)$, generically, there are jump deformations to $d_3$, $d_4$ and $d_5$, as well as smooth deformations along the family in
a neighborhood of $d_{75}(p:q)$.
\subsubsection{$(p:q)=(1:0)$}
 The algebra $d_{75}(1:0)$ has additional jump deformations to
the algebras $d_{22}\mcom d_{30}$  and the algebras $d_{45}\mcom d_{48}$.
%3, 4, 5, 22, 23, 24, 25, 26,27, 28, 29, 30,45, 46, 47, 48
\subsubsection{$(p:q)=(1:-1)$}
The algebra $d_{75}(1: -1)$ has additional jump deformations to
$d_{14}$, $d_{15}$, $d_{16}$, $d_{22}$, $d_{23}$, $d_{27}$, $d_{28}$, $d_{29}$, $d_{30}$, $d_{58}$, and  $d_{59}$.
%14, 15, 16, 22, 23, 27, 28, 29, 30, 58, 59
\subsubsection{$(p:q)=(1:1)$}
The algebra $d_{75}(1:1)$ has additional jump deformations to
$d_{22}$, $d_{23}$, $d_{27}$, $d_{28}$, $d_{29}$, $d_{30}$, $d_{37}(x: y)$ for all $(x:y)$ except $(1: 1)$ and $(0: 0)$,
$d_{38}(x: y)$ for all $(x:y)$ except $(1: 1)$ and $(0: 0)$, and $d_{73}$.
%22, 23, 27, 28, 29, 30, 37.1,37.3,37.5,38.1,38.3,38.5,73
\subsection{$(p:q)=(0:0)$}
The algebra $d_{75}(0:0)$ representing the ''generic element'' in $\P^1$, will automatically have jump deformations to every element in
the family $d_{75}(p:q)$ except itself.  In fact, it is not uncommon for the generic element to coincide with an element in some other
family, and this happens here, as $d_{75}(p:q)\sim d_{86}(1: 1)$. As a consequence, we also know that $d_{75}(0:0)$ will deform
in a neighborhood of $d_{86}(1:1)$.  Our general method of listing the algebras is organized around the principal that the order should be
preserved by deformations in the sense that elements only deform to elements whose numbers are lower.  However, the generic element in a
family cannot be expected to preserve that principle, as in some sense, it is on a higher level than members of its family.

In addition, we also expect the generic element to have jump deformations to all algebras to which any other member of its family jumps.
In fact, the algebra $d_{75}(0: 0)$ has jump deformations to
$d_1\mcom d_9$, $d_{14}$, $d_{15}$, $d_{16}$, $d_{20}\mcom d_{36}$, $d_{37}(x:y)$ $d_{38}(x:y)$ $d_{39}$, $d_{40}$,
$d_{45}\mcom d_{50}$, $d_{57}\mcom d_{60}$, $d_{65}$, $d_{68}$, $d_{73}$, $d_{74}$, $d_{75}(x:y)$, $d_{76}$, $d_{77}$, $d_{78}(x:y)$ except $(0:0)$,
$\yyy(1:1)$, $d_{80}$, $d_{82}\mcom d_{85}$, as well as deforming in a
neneighborhood of $\yyy(1: 1)$.

\subsection{Type 76}
The algebra $d_{76}=\psa{31}2+\psa{13}2+\psa{33}1+\psa{44}2$ has a 1-dimensional kernel $M=\langle e_2\rangle$. Its core
is the 3-dimensional nilpotent algebra $d_{20}(0:0)=\psa{33}1$. It is commutative.
It has jump deformations to
$d_2$, $d_6$, $d_7$, $d_8$, $d_9$, $d_{31}$, $d_{32}$, $d_{35}$, $d_{36}$, $d_{37}(1: 1)$, $d_{38}(1: 1)$, $d_{49}$, $d_{50}$, $d_{73}$, and  $d_{74}$,
as well as deforming smoothly in a neighborhood of $d_{37}(1: 1)$ and $d_{38}(1:1)$. With the exception of $d_{73}$, the algebras it jumps to are all commutative,
but the algebras $d_{37}(x:y)$ and $d_{38}(x:y)$ to which it deforms smoothly are not in general commutative.
%2, 6, 7, 8, 9, 31, 32, 35, 36, 37.2, 38.2, 49, 50, 73, 74
\subsection{Type 77}
The algebra $d_{77}=\psa{31}2+\psa{13}2+\psa{33}1+\psa{34}2$ has a 1-dimensional kernel $M=\langle e_2\rangle$. Its core is the
3-dimensional nilpotent algebra $d_{20}(0:0)=\psa{33}1$. It is isomorphic to its opposite algebra.
It has jump deformations to
$d_3$, $d_4$, $d_5$, $d_{22}\mcom d_{30}$, $d_{37}(x: y)$ except $(1: 1)$ and $(0: 0)$,
$d_{38}(x: y)$ except $(1: 1)$ and $(0: 0)$, $d_{45}\mcom d_{48}$, $d_{58}$, $d_{59}$, and $d_{73}$.
%3, 4, 5, 22, 23, 24, 25, 26, 27, 28, 29, 30,37.1,37.3,37.5,38.1, 38.3,38.5,45, 46, 47, 48, 58, 59
\subsection{Type $\mathbf{78(p:q)}$}
The family of algebras $d_{78}(p:q)$ are given by
$\psa{33}2+p\psa{41}2+q\psa{14}2+\psa{34}2+\psa{44}2$. The family is parameterized by $\P^1/2$, but the linear automorphism which transforms
$d_{78}(p:q)$ to $d_{78}(q:p)$ cannot be given in a completely generic form.  When $p=1$ and $q=0$ or $p=1$ and $q=2$, the matrix which gives the
equivalence does not fall into the same pattern as will work generically. None of the algebras in this family are commutative, owing to the
$\psa{34}2$ term in the algebra.

Generically, the kernel of this algebra is the 1-dimensional ideal $M=\langle e_2\rangle$, and its core is the trivial 3-dimensional algebra.
For the element $d_{78}(0:0)$, the dimension of the kernel is 2, and we also have that $d_{78}(0:0)\sim d_{86}(1:\gamma)$, where $\gamma$ is a
primitive 6-th root of unity.

The special cases of the parameters $(p:q)$ include not only the standard values $(1:1)$, $(1: -1)$, $(1:0)$ and $(0: 0)$, but also the values
$(1:2)$ and $(1:\gamma)$, where $\gamma$ is a primitive 6th root of unity.
\subsubsection{Generic Case}
Generically, an algebra in the family $d_{78}(p:q)$ has jump
deformations to $d_1$, $d_{37}(p:q)$ and $d_{38}(p:q)$, as well as
smooth deformations in a neighborhood of $d_{37}(p:q)$, $d_{38}(p: q)$ and $d_{78}(p:q)$.
\subsubsection{$(p:q)=(1:2)$} We mentioned something unusual in this special case already, namely that the matrices which give rise to the isomorphisms between
$d_{78}(1:2)$ and $d_{78}(2:1)$ do not fit the generic pattern. Other than this, the deformation pattern is entirely generic.
\subsubsection{$(p:q)=(1:0)$} The algebra $d_{78}(1: 0)$, in addition to the generic deformations has jump deformations to
$d_1$, $d_3$, $d_4$, $d_5$, $d_{22}\mcom d_{30}$, $d_{39}$, $d_{40}$, and $d_{45}\mcom d_{48}$.
%{1, 3, 4, 5, 22, 23, 24, 25, 26, 27, 28, 29, 30, 37(1: 0), 38(1: 0), 39, 40, 45, 46, 47, 48
\subsubsection{$(p:q)=(1:-1)$} The algebra $d_{78}(1: -1)$ has jump deformations to $d_{20}$ and $d_{21}$ in addition to the generic deformations.
\subsubsection{$(p:q)=(1:1)$} The algebra $d_{78}(1: 1)$ does not quite
follow the generic pattern in that while it does deform in a
neighborhood of $d_{37}(1:1)$
and $d_{38}(1:1)$, it does not jump to either one of them. In fact, since both $d_{37}(1: 1)$ and $d_{38}(1:1)$ are commutative, while
$d_{78}(1:1)$ is not commutative, it has no possibility of having a jump deformation to either of them. It does have a jump deformation to $d_1$, which fits the generic pattern.
\subsection{Type 79}
The algebra $\xxx=\psa{31}2+\psa{13}2+\psa{33}1$ has a 2-dimensional kernel $M=\langle e_2,e_4\rangle$. Its core is the
nontrivial 2-dimensional nilpotent algebra $d_6$. It is commutative.
It has jump deformations to
$d_2$-- $d_9$, $d_{22}$--$d_{36}$, $d_{37}(x: y)$, $d_{38}(x: y)$, both for all $(x:y)$, $d_{45}$-$d_{50}$, $d_{58}$,
$d_{59}$, $d_{65}$, $d_{68}$, $d_{73}$, $d_{74}$, $d_{76}$, and $d_{77}$.
\subsection{Type 80}
The algebra $d_{80}=\psa{33}2+\psa{43}1-\psa{34}1+\psa{44}2$ has a 2-dimensional kernel $M=\langle e_1,e_2\rangle$ and
its core is the trivial nilpotent 2-dimensional algebra.  In fact, this element is the $(1:-1)$ case of the family of
algebras $\psa{33}2+p\psa{43}1+q\psa{34}1+\psa{44}2$, which correspond to $d_{75}(p:q)$ in most cases.

The reason this element is
relegated to a different position than one would expect is that there
are two algebras which deform in a neighborhood of $d_{75}(1: -1)$. To decide which
one actually belongs in the family, we consider the following data. First, the one in the family has the smaller cohomology $H^2$. Secondly, this element
deforms to the other one, and we would expect that the one which does not deform to the other actually belongs in the family.

This algebra has jump deformations to $d_1$, $d_3$, $d_4$, $d_5$, $d_{14}$, $d_{15}$, $d_{16}$, $d_{20}\mcom d_{23}$, $d_{27}\mcom d_{30}$, $d_{39}$
$d_{40}$, $d_{57}\mcom d_{60}$, and $d_{75}(1:-1)$.
%1, 3, 4, 5, 14, 15, 16, 20, 21, 22, 23, 27, 28, 29, 30, 39, 40, 57, 58, 59, 60, 75.5
\subsection{Type 81}
The algebra $d_{81}=\psa{13}2-\psa{31}2+\psa{44}2$ has a 1-dimensional kernel $M=\langle e_2\rangle$. Its core is the
3-dimensional trivial nilpotent algebra. It is isomorphic to its opposite algebra.
It has jump deformations to
$d_1$, $d_{20}$, $d_{21}$, $d_{37}(1: -1)$, $d_{38}(1: -1)$, $d_{40}$,
$d_{51}$, $d_{52}$, $d_{78}(1: -1)$ and deforms in neighborhoods of
the algebras $d_{37}(1:-1)$. $d_{38}(1:-1)$ and $d_{78}(1:-1)$.
\subsection{Type 82}
The algebra $d_{81}=\psa{33}2+\psa{43}1+\psa{34}1+\psa{44}2$ has a 2-dimensional kernel $M=\langle e_1,e_2\rangle$. Its core is the
2-dimensional trivial nilpotent algebra. It is isomorphic to its opposite algebra.
It has jump deformations to
$d_3$, $d_4$, $d_5$, $d_{22}$, $d_{23}$,
$d_{27}\mcom d_{30}$, $d_{37}(x:y)$ except $(1:1)$ and $(0:0)$,
$d_{38}(x:y)$ except $(1:1)$ and $(0:0)$, $d_{73}$
$d_{75}(1:1)$ and deforms in a neighborhood of $d_{75}(1:1)$.
%3, 4, 5, 22, 23, 27, 28, 29, 30,37.1,37.3,37.5,38.1,38.3,38.5, 73
\subsection{Type $\mathbf{83(p:q)}$}
The family of algebras $\yyy(p:q)$ are given by 
$p\psa{42}1+q\psa{24}1+\psa{44}1+\psa{44}3$. The family is parameterized by $\P^1$, and does not have an action of $\Z_2$, so that $\yyy(p:q)$ is not, in general, isomorphic to $\yyy(q:p)$.
Generically, the kernel of this algebra is the 2-dimensional ideal $M=\langle e_1,e_3\rangle$, and its core is the trivial 2-dimensional algebra.
For the element $\yyy(0:0)$, the dimension of the kernel is 3, and we also have that $\yyy(0:0)\sim d_{86}(0:0)$.

The special cases of the parameters $(p:q)$ include not only the standard values $(1:1)$, $(1: -1)$, $(1:0)$ and $(0: 0)$, but also $(0:1)$ as $\yyy(1: 0)$ is not isomorphic to $\yyy(1:0)$, so must be treated separately.
\subsubsection{Generic Case}
Generically, an algebra in the family $\yyy(p:q)$ has jump deformations to
$d_3$, $d_4$, $d_5$, $d_{14}$, $d_{15}$, $d_{16}$, $d_{22}\mcom d_{30}$ $d_{37}(x:y)$ except $(1: 1)$ and $(0: 0)$, $d_{38}(x:y)$, except $(1:1)$ and $(0:0)$,
 $d_{45}\mcom d_{48}$, $d_{58}$, $d_{59}$, $d_{73}$, $d_{75}(x: y)$
 except $(0:0)$, and $d_{77}$ as well as smooth deformations in a
 neighborhood of  $\yyy(p:q)$.
%3, 4, 5, 14, 15, 16, seq(i,i=22..30), 37.1,37.3,37.5,38.1,38.3,38.5,45,46,47, 48, 58, 59, 73, 75.1,75.2,75.3,75.5, 77

\subsubsection{$(p:q)=(1:0)$} The algebra $\yyy(1: 0)$, in addition to the generic deformations has jump deformations to
$d_{11}$, $d_{12}$, $d_{42}$, $d_{44}$, $d_{67}$, and $d_{69}$.
%3, 4, 5, 11, 12, 14, 15, 16, 22, 23, 24, 25, 26, 27, 28, 29, 30, 37.1,37.3,37.5,38.1,38.3,38.5, 42, 44, 45, 46, 47, 48, 58, 59, 67, 69, 73, 75.1,75.2,75.3,75.5,77
\subsection{$(p:q)=0:1)$} The algebra $\yyy(0:1)$, in addition to the generic deformations has jump deformations to
$d_{10}$, $d_{13}$, $d_{41}$, $d_{43}$, $d_{66}$, and $d_{70}$. Note that these are precisely the opposite algebras to the extra deformations of the algebra
$\yyy(1:0)$.
%3, 4, 5, 10, 13, 14, 15, 16, 22, 23, 24, 25, 26, 27, 28, 29, 30, 37.1,37.3,37.5,38.1,38.3,38.5, 41, 43, 45, 46, 47, 48, 58, 59, 66, 70, 73, 75.1,75.2,75.3,75.5, 77
\subsection{$(p:q)=(1:-1)$} The algebra $\yyy(1: -1)$ has jump deformations to $d_{1}$,
$d_{20}$, $d_{21}$, $d_{39}$, $d_{40}$, $d_{51}$, $d_{52}$, $d_{60}$ and $d_{80}$, in addition to the generic deformations.
%1, 3, 4, 5,14, 15, 16, 20, 21, 22,23,24, 25, 26, 27, 28,29,30, 37.1,37.3,37.5, 38.1,38.3,38.5, 39, 40, 45,46,47,48,51, 52, 57, 58,59,60,73, 75.1,75.2,75.3,75.5, 77, 80
\subsection{$(p:q)=(1:-1)$} In addition to the generic deformations, the algebra $\yyy(1: 1)$ has jump deformations to
$d_2$, $d_6$, $d_7$, $d_8$, $d_9$, $d_{31}\mcom d_{36}$, $d_{37}(1:1)$, $d_{38}(1:1)$, $d_{49}$, $d_{50}$, $d_{65}$, $d_{68}$, $d_{74}$, and $d_{82}$.
Note that except for $d_{65}$ all of these algebras are commutative.
%2, 3, 4, 5, 6, 7, 8, 9, 14, 15, 16, 22, 23, 24,25, 26, 27,28, 29,30,31, 32,33, 34, 35, 36, 37.1,37.2,37.3,37.5,38.1, 38.2,38.3,38.5,38.6, 45,46,47,48,49, 50, 58, 59, 65, 68, 73,74, 75.1,75.2,75.3,75.5, 77, 82

\subsection{Type 84}
The algebra $d_{84}=\psa{13}2-\psa{31}2+\psa{14}2$ has a 1-dimensional kernel $M=\langle e_2 \rangle$. Its core is the trivial
3-dimensional nilpotent algebra.
It has jump deformations to
$d_1$, $d_3$, $d_4$, $d_5$, $d_{14}$, $d_{15}$, $d_{16}$,$d_{20}$--$d_{30}$,
$d_{37}(x: y)$ except $(1: 1)$ and $(0: 0)$, $d_{38}(x: y)$ except $(1: 1)$ and $(0: 0)$,
$d_{39}$, $d_{40}$, $d_{45}$--$d_{48}$, $d_{57}$--$d_{60}$, $d_{73}$, $d_{77}$, $d_{78}(x:y)$ except $(0:0)$.
\subsection{Type 85} The algebra $d_{85}=\psa{31}2+\psa{13}2+\psa{33}2+\psa{44}2$ has a 1-dimensional kernel $M=\langle e_2 \rangle$. Its core is the trivial
3-dimensional nilpotent algebra. It has jump deformations to $d_1\mcom d_{20}$, $d_{22}$, $d_{23}$, $d_{27}\mcom d_{32}$, $d_{35}$, $d_{36}$,
$d_{37}(x:y)$ except $(0:0)$, $d_{38}(x:y)$ except $(0:0)$, $d_{73}$.
$d_{74}$, $d_{76}$, $d_{78}(1:1)$, and deforms in a neighborhood of $d_{78}(1:1)$.
%seq(i,i=1..9), 20, 22, 23, seq(i,i=27..32), 35, 36, 37.1,37.2,37.3,37.5,38.1,38.2,38.3,38.5, 49, 50, 73, 74, 76, 78.2
\subsection{Type $\mathbf{86(p:q)}$}
The family of algebras $d_{86}(p:q)=\psa{42}1p+\psa{24}1q+\psa{44}1$ is parameterized by $\P^1/\Sigma_2$. The algebras are not commutative,
except for $d_{86}(1:1)$ and $d_{86}(0:0)$. Generically, the elements in this family have a 2-dimensional kernel $M=\langle e_1,e_3\rangle$,
and their core is the trivial 2-dimensional algebra, except that $d_{86}(0:0)$ has 3-dimensional kernel and 1-dimensional core.
\subsubsection{Generic Case}
Generically, $d_{86}(p:q)$ has jump deformations to $d_1$, $d_3$, $d_4$, $d_5$, $d_{14}$, $d_{15}$, $d_{16}$, $d_{20}$--$d_{30}$, $d_{37}(x: y)$
except $(1: 1)$ and $(0: 0)$, $d_{38}(x: y)$ except $(1: 1)$ and $(0: 0)$, $d_{39}$, $d_{40}$ $d_{45}$--$d_{48}$, $d_{57}$--$d_{60}$, $d_{73}$,
$d_{75}(x: y)$ except $(0: 0)$, $d_{77}$, $d_{78}(x: y)$ except $(0:
0)$, $\yyy(p: q)$, $d_{80}$, and $d_{84}$, as well as deforming
smoothly in a neighborhood of $\yyy(p: q)$ and $d_{86}(p: q)$.
%1, 3, 4, 5, 14, 15, 16, seq(i,i=20..30), 37.1, 37.3, 37.5, 38.1, 38.3, 38.5, 39, 40, 45, 46, 47, 48, 57, 58, 59, 60, 73,  75.1, 75.2, 75.3, 75.5, 77, 78.1, 78.2, 78.3, 78.5, 79.1, 80, 84
\subsubsection{$(p:q)=(1:0)$}
The algebra $d_{86}(1:0)$ has additional jump deformations to $d_{10}\mcom d_{14}$, $d_{41}$, $d_{42}$, $d_{43}$, $d_{61}$, $d_{62}$, $d_{63}$, $d_{77}$ and $d_{84}$.
%1, 3, 4, 5, 14, 15, 16, seq(i,i=20..30), 37.1, 37.3, 37.5, 38.1, 38.3, 38.5, 39, 40, 45, 46, 47, 48, 57, 58, 59, 60, 73,  75.1, 75.2, 75.3, 75.5, 77, 78.1, 78.2, 78.3, 78.5, 79.1, 80, 84
\subsubsection{$(p:q)=(1:1)$}
The algebra $d_{86}(1:1)$ has additional jump deformations to $d_2$, $d_6$, $d_7$, $d_8$, $d_9$, $d_{31}\mcom d_{36}$, $d_{49}$, $d_{50}$, $d_{65}$, $d_{68}$,
$d_{74}$, $d_{76}$, $d_{82}$ and $\xxx$. With the exception of $d_{65}$, all of these additional jump deformations are to commutative algebras, to which it is not possible for generic elements of the family to deform.
\subsubsection{$(p:q)=(1:-1)$}
The algebra $d_{86}(1:1)$ has additional jump deformations to $d_{51}$, $d_{52}$ and $d_{81}$.
%1, 3, 4, 5, 14, 15, 16, 20, 21, 22, 23, 24, 25, 26, 27, 28, 29, 30,37.1,37.3,37.5,38.1,38.3,38.5, 39, 40, 45, 46, 47, 48, 51, 52, 57, 58, 59, 60, 73,75.1,75.2,75.3,75.5,77, 78.1,78.2,78.3,78.5, 79.5, 80, 81, 84
\subsection{Type 87} The algebra $d_{87}=-\psa{31}2+\psa{13}2$ has a 2-dimensional kernel $M=\langle e_2,e_4 \rangle$. Its core is the trivial
2-dimensional nilpotent algebra. It has jump deformations to $d_1$, $d_3$, $d_4$, $d_5$,$d_{14}$, $d_{15}$, $d_{16}$, $d_{20}\mcom d_{30}$,
$d_{37}(x:y)$ except $(1:1)$ and $(0:0)$, $d_{38}(x:y)$ except $(1:1)$ and $(0:0)$, $d_{39}$, $d_{40}$, $d_{45}\mcom d_{48}$, $d_{51}\mcom d_{54}$,
$d_{57}\mcom d_{60}$, $d_{73}$, $d_{75}(x:y)$ except $(0:0)$, $d_{77}$,
$d_{78}(x:y)$ except $(0:0)$, $\yyy(1:-1)$, $d_{80}$, $d_{81}$ $d_{84}$
and $d_{86}(1:-1)$ and deforms in a neighborhood of $\yyy(1: -1)$ and $d_{86}(1:-1)$.
%1, 3, 4, 5, 14, 15, 16, 20, 21, 22, 23, 24, 25, 26, 27, 28, 29, 30, 37.1,37.3,37.5,38.1,38.3,38.5,39, 40, 45, 46, 47, 48, 51, 52, 53, 54, 57, 58, 59, 60, 73, 75.1,75.2,75.3,75.5,77,78.1,78.2,78.3,78.5,79.5, 80, 81, 84, 86.5]
\subsection{How we constructed the algebras}
In section 1, we discussed that fact that any nilpotent algebra has a unique kernel, consisting of all elements whose product with any element is zero, and core, consisting of the quotient of the algebra by its kernel, so that, in principle, it should be easier to construct all algebras of a given dimension by considering extensions of algebras by ideals in this manner.

However, our motivation for the constructions is not simply to give a list of all algebras, but to give a decomposition of the moduli space into strata which are dictated by deformation theory. A single stratum may include algebras whose kernels do not all have the same dimension, so this fact alone dictates the need for another method of identifying the strata.  Moreover, our constructions have revealed that the strata seem to consist of projective orbifolds, and it is necessary to identify these orbifolds in some natural way.  This has led us to the following method of
construction.

We first need to know all $(n-1)$-dimensional nilpotent algebras, and then we use the fact that there is always a codimension 1 ideal in any nilpotent algebra to realize all algebras as extensions of the trivial 1-dimensional algebra by an $(n-1)$-dimensional nilpotent algebra. In the case of dimension 4 complex nilpotent algebras, there are 3 nontrivial nilpotent 3-dimensional algebras, as well as the trivial nilpotent algebra, which come into play.  Here we will consider how our algebras arise by this method of construction.
\subsection{Extensions by the nilpotent algebra $d_{19}=\psa{31}2+\psa{13}2+\psa{33}1$} The algebras $d_{73}$, $d_{76}$, $d_{77}$, and $\xxx$ all arise
as extensions of the trivial algebra by the 3-dimensional algebra $d_{19}$. From  compatibility condition $[\mu,\lambda]=0$ and the Maurer Cartan equation, we can reduce  $\lambda$ to the form $\lambda=\psa{34}2x$, where $x$ can be chosen  to be either 1 or 0, $\psi$ can be taken to be 0, and $\tau$ is of the form
$\tau=\psa{44}2c$, where again, $c$ is either 0 or 1. This gives 4 distinct possibilities, which are the four algebras listed above.

\subsection{Extensions by the nilpotent algebra $d_{20}(p:q)=\psa{13}2p+\psa{31}2q+\psa{33}2$} The algebras
$d_{73}$, $d_{74}$, $d_{75}(x: y)$, $d_{76}$, $d_{78}(x: y)$, $\yyy(1: 1)$, $\yyy(0: 0)$, $d_{80}$,
$d_{81}$, $d_{82}$, $d_{84}$, $d_{85}$, and $d_{86}(x: y)$ arise as extensions of the trivial algebra by the nilpotent algebra $d_{20}(p:q)$.
In this case, the compatibility condition gives rise to more than one solution, and these solutions depend on the variable $(p:q)$. For the special values
$(1:0)$ and $(0:0)$, there are nongeneric solutions in addition to the generic case.

Let us examine the generic case first. In this case, we can reduce
$\lambda$ to the form $\lambda=\psa{43}2r+\psa{14}2s$, where $(r:s)$ is
a projective coordinate. However, after further analysis, we obtain that the algebras which arise in this fashion are
$d_{81}$, $d_{84}$, $d_{85}$ and $d_{86}(p:q)$. Thus the projective coordinates $(r:s)$ don't arise in the final description of the algebra.

When $(p:q)=(1:0)$, we obtain the same format for $\lambda$ and $\tau$, and the algebras we obtain
are $d_{78}(x:y)$ except $(0:0)$, $d_{84}$ and $d_{86}(1: 0)$.

Finally, when $(p:q)=(0:0)$, we obtain 4 distinct solutions for the format of $\lambda$. For example, in one of the patterns, we have
$\lambda=\psa{41}2p+\psa{14}2q +\psa{34}2u$, where $(p:q)$ is a projective coordinate, and the choice of $u$ can be reduced to $u=1$ or $u=0$. Thus
we obtain a new projective coordinate, and this coordinate is retained in some of the solutions, so we obtain new projective families of algebras which
don't arise from the old projective families.  After some analysis, we obtain that algebras arising in this case are
$d_{73}$, $d_{74}$, $d_{75}(x:y)$, $d_{76}$, $d_{78}(x: y)$, $\yyy(1:1)$, $\yyy(0: 0)$, $d_{80}$, $d_{82}$, and  $d_{85}$.
%73, 74, 75(x: y), 76, 78(x: y), 79(1: 1), 79(0: 0)=86(0: 0), 80, 82, 85, 86(1: 1), 86(1: alpha), 86(0: 0)
\subsection{Extensions by the nilpotent algebra $d_{21}=\psa{13}2-\psa{31}2$}
The algebras $d_{78}(1: -1)$, $d_{81}$, $d_{84}$, and  $d_{87}$ all arise as extensions of the nilpotent 3-dimensional algebra $d_{21}$.
\subsection{Extensions by the trivial 3-dimensional nilpotent algebra} The algebras $d_{77}$, $\yyy(x: y)$, $d_{80}$, $\xxx$, $d_{84}$, $d_{86}(x: y)$,
and $d_{87}$ arise as extensions of the trivial 1-dimensional algebra by the trivial 3-dimensional nilpotent algebra.

From the description above, it is clear that many of these algebras arise as extensions of the trivial algebra by different 3-dimensional nilpotent algebras.

\section{Commutative Algebras}

There are 20 distinct non-nilpotent commutative algebras, of which 9 are unital. Every commutative algebra is a direct sum of
algebras which are ideals in quotients of polynomial algebras. Every finite dimensional unital commutative algebra is a quotient of a polynomial algebra, while every
finite dimensional nonunital algebra is an ideal in such an algebra. 
In \cite{fp14} we gave a table showing the non-nilpotent commutative algebras,
which we will not reproduce here.

Nilpotent commutative algebras were
classified by Hazlett \cite{hazl}, and also given in \cite{mazz2}. There are 8 nontrivial nilpotent commutative algebras.

We note that commutative algebras may deform into noncommutative algebras, but noncommutative algebras never deform into
a commutative algebras.  The fact that commutative algebras have noncommutative deformations plays an important role in physics,
and deformation quantization describes a certain type of deformation of a commutative algebra into a noncommutative one.

\begin{table}[h,t]
\begin{center}
\begin{tabular}{ll}
Algebra&Structure\\ \hline \\
$d_{74}$&$x\C[x]/(x^5)$\\
$d_{75}(1:1)$&$(x,y)\le \C[x,y]/(x^2-y^2,yx^2,xy^2) $\\
$d_{75}(0:0)=d_{86}(1:1)$&$\C_0\oplus(x,y)\le\C_0\oplus\C[x,y]/(x^2-y^2,xy)$\\
$d_{76}$&$(x,y)\le\C[x,y]/(y^3-x^2,xy)$\\
$\xxx$&$\C_0\oplus x\C[x]/(x^4)$\\
$\yyy(1:1)$&$(x,y)\le\C[x,y]/(y^2,x^2y,x^3)$\\
$\yyy(0:0)=d_{86}(0:0)$&$\C_0^2\oplus x\C[x]/(x^3)$\\
$d_{85}$&$(x,y,z)\le\C[x,y,z]/(x^2-y^2,y^2-yz,xy,xz,z^2)$\\
%$d_{86}(1:1)$&$\C_0\oplus(x,y)\le\C_0\oplus\C[x,y]/(x^2-xy,y^2)$\\
$d_0$&$\C_0^4$\\
\end{tabular}
\end{center}
\label{Commutative Nilpotent}
\caption{Nilpotent 4-dimensional commutative algebras}
\end{table}

\section{Levels of algebras}
We give a table showing the levels of each algebra. For completeness,
we include the levels of the non-nilpotent
algebras, so that the reader can see how the nilpotent algebras fit into the general picture.
To define the level,
we say that a rigid algebra has level 1, an algebra which has only
jump deformations to an algebra on level one has level two and so on. To be on level $k+1$,
an algebra must have a jump deformation to an algebra on level $k$, but no jump deformations
to algebras on a level higher than $k$. For families, if one algebra in the family has a jump
to an element on level $k$, then we place the entire family on at least level $k+1$. Thus, even
though generically, elements of the family $d_{37}(p:q)$ deform only to members of the same
family, there is an element in the family which has a jump to an element on level 4. For
the generic element in a family, we consider it to be on a higher level
than the other elements, 
because it has jump deformations to the other elements in its family.

\begin{table}[h,t]
\begin{center}
\begin{tabular}{ll}
Level&Algebras\\ \hline \\
$1$&1,2,3,4,5,10,11,12,13,14,15,16,17,18,19,20,39,53,54,55,56\\
$2$&6,7,21,22,23,24,25,26,27,28,29,30,40,41,42,43,44,57\\
$3$&8,9,31,32,45,46,47,48,58,59,60,61,62,63,64\\
$4$&33,34,35,36,49,50,66,67,69,70\\
$5$&$37(p:q)$, $38(p:q)$,65,68,74\\
$6$&$37(0:0)$,$38(0:0)$,51,52,73,${75}(p:q)$\\
$7$&71,72,76,77,$78(p:q)$,80,82,85\\
$8$&79,81,$83(p:q)$,84\\
$9$&$75(0:0)$,$78(0:0)$,$86(p:q)$\\
$10$&$83(0:0)$,$86(0:0)$,87\\
\\ \hline
\end{tabular}
\end{center}
\label{Levels}
\caption{The levels of the algebras}
\end{table}
\section{Analysis}
The reader may note that every element in the family $d_{86}(p:q)$ has
a jump deformation to every element in the families $d_{75}(x:y)$ and
$d_{78}(x:y)$. As a consequence, it is difficult to say which element
in the family $d_{86}(p:q)$ really should be the $(0:0)$ element in
those families. Of course, for the choice of algebras representing
$d_{75}(p:q)$, the corresponding member of the family is precisely the
element $d_{86}(1:1)$, but this can be an artifact of the form of the
algebra structure.

To illustrate this point, we consider the following family of algebras, given by
$d(p:q)=\psa{13}2+\psa{33}2+p\psa{41}2+q\psa{34}2$. It is easily checked that $d(q:p)\sim d(p:q)$, and that the family is given by projective coordinates. In fact, $d(p:q)\sim d_{78}(x:y)$ where $\frac qp=\frac{xy}{x^2+xy+y^2}$, so in general, the family $d(p:q)$ and $d_{78}(x:y)$ determine the same collection of algebras.  However,
$d(0:0)$ is isomorphic to $d_{86}(1:0)$, not $d_{86}(1:\alpha)$ (where
$\alpha$ is a primitive 6th root of unity), so we see that there may be
no natural way to identify the generic element in a projective family
as a specific algebra.

What is important is that there always is an algebra corresponding to the 'so called' generic element of the $\P^1$ which parameterizes a family.

There are two additional features of this moduli space that did not
arise in the lower dimensional moduli spaces of ordinary complex
associative algebras. The first is that we were unsuccessful it
describing the family $d_{75}(p:q)$ with a single family of algebras
with parameters $p$ and $q$. (A similar difficulty arose with the
moduli space of $2|2$-dimensional \zt-graded associative
algebras,\cite{fop1}, so this pattern has been seen before, just not in
the context of ordinary associative algebras.)

The second feature we saw was that for the algebras $d_{73}$ and
$d_{75}(1:1)$, there were deformations in a neighborhood of elements of
another family, without there being a jump to the point in whose
neighborhood the smooth deformation occurred. This was a bit surprising to us, because we had not observed this behavior before.

The authors and collaborators have written a series of articles describing moduli spaces of Lie, associative, \linf, and \ainf\ algebras.  It first became
clear when studying the 3-dimensional complex Lie algebras that there was a stratification of the moduli space by projective orbifolds, each of the form
$\P^n$, possibly with an action of $\Sigma_{n+1}$, given by interchanging coordinates. Our example of 4-dimensional complex associative algebras fits this
pattern nicely.  The authors have been conjecturing that this type of stratification by projective orbifolds ought to hold in general, but as of yet, do not have
the tools to prove this.

%\bibliographystyle{amsplain}
%\bibliography{global}
\providecommand{\bysame}{\leavevmode\hbox to3em{\hrulefill}\thinspace}
\providecommand{\MR}{\relax\ifhmode\unskip\space\fi MR }
% \MRhref is called by the amsart/book/proc definition of \MR.
\providecommand{\MRhref}[2]{%
  \href{http://www.ams.org/mathscinet-getitem?mr=#1}{#2}
}
\providecommand{\href}[2]{#2}

\end{document}